\documentclass{elsart}
\usepackage{amssymb}

\newcommand{\To}{\longrightarrow}

\newcommand{\RR}{\mathbb{R}}
\newcommand{\Rd}{{\mathbb{R}}^d}
\newcommand{\NN}{\mathbb{N}}
\newcommand{\ZZ}{\mathbb{Z}}
\newcommand{\QQ}{\mathbb{Q}}

\newcommand{\EA}{\mathcal{E}_A}
\newcommand{\EE}{\mathcal{E}}

\newcommand{\uu}{\mathbf{u}}
\newcommand{\vv}{\mathbf{v}}

\newcommand{\xx}{\mathbf{x}}
\newcommand{\yy}{\mathbf{y}}

\newcommand{\sn}{\{\mathbf{0}\}}

\newcommand{\comment}[1]{}

\newcounter{res}
\newcounter{rev}
\begin{document}
\begin{frontmatter}

\title{Equivalence of $A$-Approximate Continuity for Self-Adjoint
Expansive Linear Maps}

\author[Hung]{Sz. Gy. R\'ev\'esz\thanksref{hu}},
\ead{revesz@renyi.hu} 
\author[Spa]{A. San Antol\'{\i}n\corauthref{cor}\thanksref{sp}}
\ead{angel.sanantolin@uam.es}
\address[Hung]{A. R\'enyi Institute of Mathematics,
Hungarian Academy of Sciences, Budapest, P.O.B. 127, 1364 Hungary}
\address[Spa]{Departamento de Matem\'aticas,
Universidad Aut\'onoma de Madrid, 28049 Madrid, Spain}
\corauth[cor]{Corresponding author. Telf: 0034 91 497 5253, fax:
0034 91 497 4889}
\thanks[hu]{Present address: Institut Henri Poincar\'e,
11 rue Pierre et Marie Curie, 75005 Paris, France}
\footnotetext{Supported in part in the framework of the
Hungarian-Spanish Scientific and Technological Governmental
Cooperation, Project \# E-38/04.} \footnotetext{The author was
supported in part by the Hungarian National Foundation for
Scientific Research, Project \# T-049301.} \footnotetext{This work
was accomplished during the author's stay in Paris under his Marie
Curie fellowship, contract \# MEIF-CT-2005-022927.}
\thanks[sp]{Partially supported by  \#
MTM 2004-00678}  \footnotetext{Supported in part in the framework of
the Spanish-Hungarian Scientific and Technological Governmental
Cooperation, Project \# HH 2004-0002.}

\begin{abstract}  Let $A:\RR^{d}\longrightarrow \RR^{d}$, $d
\geq 1$, be an expansive linear map. The notion of $A$-approximate
continuity was recently used to give a characterization of scaling
functions in a multiresolution analysis (MRA). The definition of
$A$-approximate continuity at a point $\mathbf{x}$ -- or,
equivalently, the definition of the family of sets having
$\mathbf{x}$ as point of $A$-density -- depend on the expansive
linear map $A$. The aim of the present paper is to characterize
those self-adjoint expansive linear maps $A_1, A_2: \Rd \to \Rd$ for
which the respective concepts of $A_{\mu}$-approximate continuity
($\mu=1,2$) coincide. These we apply to analyze the equivalence
among dilation matrices for a construction of systems of MRA. In
particular, we give a full description for the equivalence class of
the dyadic dilation matrix among all self-adjoint expansive maps. If
the so-called ``four exponentials conjecture'' of algebraic number
theory holds true, then a similar full description follows even for
general self-adjoint expansive linear maps, too.
\end{abstract}
\end{frontmatter}

\section{Introduction}

A multiresolution analysis (MRA) is a general method introduced by
Mallat \cite{M:89} and Meyer \cite{Me:90} for constructing
wavelets. On ${\RR}^{d} \ (d\geq 1)$ equipped with the Euclidean
norm $\|\cdot\|$, an MRA means a sequence of subspaces $V_j$,
$j\in\ \mathbb{Z}$ of the Hilbert space $L^2(\RR^d)$ that
satisfies the following conditions:
\begin{enumerate}
\item[($\mbox{i}$)] \begin{description} $\forall j\in \mathbb{Z},\qquad V_j\subset
V_{j+1}$; \end{description} \item[($\mbox{ii}$)]
\begin{description} $\forall j\in \mathbb{Z}, \qquad f({\bf x})\in
V_j \Leftrightarrow f(2{\bf x}) \in V_{j+1}$; \end{description}
\item[($\mbox{iii}$)]\begin{description} $W=\overline{\bigcup_{j\in\mathbb{Z}} V_j} =
L^2(\RR^d)$; \end{description} \item[($\mbox{iv}$)]
\begin{description}there exists a \emph{scaling function} $\phi\in
V_0$, such that $\{ \phi ({\bf x} - {\bf k}) \}_{{\bf
k}\in\mathbb{Z}^d}$ is an orthonormal basis for $V_0$.
\end{description}
\end{enumerate}

We could consider MRA in a  general context, where instead of the
dyadic dilation one considers a fixed linear map $A:
\RR^d\rightarrow \RR^d$ such that $A$ is an expansive map, i.e.
all (complex) eigenvalues have absolute value greater than 1, and
\begin{equation}\label{uno}
A(\mathbb{Z}^d) \subset \mathbb{Z}^d,
\end{equation}
i.e., the corresponding matrix of $A$ with respect to the
canonical basis has every entries belonging to $\mathbb{Z}$. Given
such a linear map $A$ one defines an $A-$MRA as a sequence of
subspaces $V_j$, $j\in\mathbb{Z}$ of the Hilbert space
$L^2(\RR^d)$ (see \cite{Ma:92}, \cite{GM:92}, \cite{S:94},
\cite{W:97}) that satisfies the conditions
$(\mbox{i}),(\mbox{iii}),(\mbox{iv})$ and
\begin{enumerate}
\item[($\mbox{ii}_1$)] \begin{description} $\forall j\in\mathbb{Z}, \qquad f({\bf
x})\in V_j \Leftrightarrow f(A{\bf x}) \in V_{j+1}.$
\end{description}
\end{enumerate}
A characterization of scaling functions in a multiresolution
analysis in a general context was given in \cite{CKS:05}, where
the notion of $A$-approximate continuity is introduced as a
generalization of the notion of approximate continuity.

In this work $|G |_{d}$ denotes the $d$-dimensional Lebesgue measure
of the set $G\subset\Rd$, and $B_{r}:= \{ \mathbf{x} \in \RR^d~:~
\parallel \mathbf{x} \parallel < r \}$ stands for the ball of
radius $r$ with the center in the origin. Also, we write
$F+\mathbf{x_0} =  \{\mathbf{y}+\mathbf{x_0} ~:~ \mathbf{y} \in F\}$
for any $F \subset \RR^d$, $\mathbf{x_0} \in \RR^{d}$.

\begin{defn}  \label{Adefn1} Let an expansive linear map
$A:\RR^{d} \longrightarrow \RR^{d}$ be given. It is said that
$\mathbf{x}_{0} \in \RR^{d}$ is a point of $A-$density for a
measurable set $E \subset \RR^d,$ $|E|_d>0$ if for all $r>0$,
\begin{equation}
\lim_{j\longrightarrow \infty}\frac{| E\bigcap
(A^{-j}B_{r}+\mathbf{x}_{0})|_{d}}{|
A^{-j}B_{r}+\mathbf{x}_{0}|_{d}} =1. \label{AD1}
\end{equation}
Given an expansive linear map $A:\RR^{d}\longrightarrow \RR^{d}$,
and given $ \mathbf{x_{0}} \in \RR^{d}$, we define the
\emph{family of $A$-dense sets at $\mathbf{x_{0}} $} as
\[
\mathcal{E}_{A}(\mathbf{x_{0}}) = \{E \subset \RR^{d} \textrm{
measurable set : $\mathbf{x_{0}}$ is a point of $A-$density for
$E$}\}.
\]
Furthermore, we will write $\mathcal{E}_{A}$ when $\mathbf{x_{0}}$
is the origin.
\end{defn}

\begin{defn} \label{A-contiaprox}
Let $A: \RR^d \to \RR^d$ be an expansive linear map and let $f:
\mathbb{R}^n \longrightarrow \mathbb{C}$ be a measurable function.
It is said that $\mathbf{x}_{0}\in \mathbb{R}^n$ is a point of
$A$-approximate continuity of the function  $f$ if there exists a
measurable set $E \subset \mathbb{R}^{n}$, $\mid E \mid_{n}>0$, such
that $\mathbf{x}_{0}$ is a point of  $A$-density for the set $E$ and
\begin{equation}\label{CCcuatro} \lim_{{\bf x} \to \mathbf{x}_0, {\bf
x}\in E} f({\bf x}) = f(\mathbf{x}_0).
\end{equation}
\end{defn}

The relation between the behavior of the Fourier transform
$\widehat{\phi}$ of the scaling function $\phi$ in the neighborhood
of the origin and the condition (iii) is described in the following
theorem of \cite{CKS:05}.

{\bf Theorem A.} Let $V_j$ be a sequence of closed subspaces in
$L^2(\RR^d)$ satisfying the conditions $(\mbox{\rm i})$, $(\mbox{\rm
ii}_1)$ and $(\mbox{\rm iv})$. Then the following conditions are
equivalent:
\begin{description}
\item[${\bf (A)}$]
$W=\overline{\bigcup_{j\in\mathbb{Z}}V_j}=L^2(\RR^d)$;
\item[${\bf (B)}$]  Setting $|\widehat{\phi}({\bf 0})| = 1$,
the origin is a point of $A^*$-approximate continuity of the
function $|\widehat{\phi}|$.
\end{description}

As it was observed in \cite[Remark 5, p. 1016]{CKS:05}, that the
definition of points of $A$-approximate continuity depends of the
expansive linear map $A$.

The aim of the present paper is to study the following problem:

{\bf Problem 1:} Characterize those expansive linear maps $A_1, A_2:
\Rd \to \Rd$ for which the concept of $A_1$-approximate continuity
coincides with the concept of $A_2$-approximate continuity.

{\bf Remark 1.} 
>From the definition of point of $A$-approximate continuity of a
measurable function on $\mathbb{R}^d$, it is easy to see that
given $\mathbf{x_0} \in \mathbb{R}^d$ and given a measurable set
$E \subset \mathbb{R}^d$, the point $\mathbf{x_0}$ is a point of
$A$-approximate continuity for the function
\[
 f(\mathbf{x})=
\chi_{E\cup  \{ \mathbf{x_{0}}\}}(\mathbf{x})= \left \{
\begin{array}{ll}
1 & \textrm{if $\mathbf{x} \in E \bigcup \{ \mathbf{x_{0}} \}$} \\
0 & \textrm{if $\mathbf{x} \notin E$}
\end{array}  \right.
\]
if and only if $E \in \mathcal{E}(\mathbf{x_{0}})$. Therefore, it
suffices to study the notion of $A$-density, and once
$\mathcal{E}_{A_1}(\mathbf{x_0})=\mathcal{E}_{A_2}(\mathbf{x_0})$,
also the notions of $A_1$-approximate continuity and
$A_2$-approximate continuity coincide.

Moreover, clearly, $E \in \mathcal{E}_{A}$ if and only if $E +
\mathbf{x_0} \in \mathcal{E}_{A}(\mathbf{x_0})$.

Thus, we can simplify Problem 1 in the following way.

{\bf Problem 1':} Describe under what conditions on two expansive
linear maps $A_1, A_2: \mathbb{R}^d \To \mathbb{R}^d$, we have
that $\mathcal{E}_{A_1}=\mathcal{E}_{A_2}$.

In Corollary \ref{cormain} we solve the problem for
\emph{expansive self-adjoint linear maps} on $\mathbb{R}^d$,
without the extra condition (\ref{uno}). In the last section we
discuss the additional, essentially number theoretical
restrictions, brought into play by condition (\ref{uno}).

Characterization of expansive matrices satisfying (\ref{uno}) have
been studied by several authors.

In \cite{LW:95} a complete classification for expanding $2 \times
2$-matrices satisfying (\ref{uno}) and $\mid \det M \mid =2$ is
given. Their result is the following.

Call two integer matrices $A$ and $M$ integrally similar if there
exists an integer unimodular matrix $C$ such that $C^{-1}AC= M$.
Now, denote

 $A_1= \left(\begin{array}{rr} 0 & 2 \\ 1 & 0
\end{array}\right)\;\;\;\;$
 $A_2= \left(\begin{array}{rr} 0 & 2 \\ -1 & 0
\end{array}\right)$ \\

$A_3= \left(\begin{array}{rr} 1 & 1 \\ -1 & 1
\end{array}\right)\;\;\;\;$
$A_4= \left(\begin{array}{rr} 0 & 2 \\ -1 & 1
\end{array}\right).$

{\bf Lemma B.} {\em Let $M$ be an expanding  $2 \times 2$-matrix
satisfying (\ref{uno}). If   $\det M =-2$ then $M$ is integrally
similar to $A_1$. If $\det M =2$ then $M$ is integrally similar to
$A_2$, $\pm A_3$, $\pm A_4$. }

On the other side, a complete characterization for expanding $2
\times 2$-matrices satisfying (\ref{uno}) and
\begin{equation}\label{4.1}
M^l=nI \qquad \textrm{for some} \qquad l,n \in \mathbb{N}
\end{equation}
is given in \cite{DLN:97}. Their answer is given in the following
theorem where they do not write the trivial case that $M$ is a
diagonal matrix.

{\bf Theorem C.} {\em Given $l,n \in \mathbb{N}$, an expanding $2
\times 2$-matrix which satisfy (\ref{uno}) and  $(\ref{4.1})$
exists if and only if there exist two numbers $\lambda_1 \neq
\lambda_2$ whose sum and product are integral and satisfy
$\lambda_1^l= \lambda_2^l = n$. Furthermore, then $\lambda_1$ and
$\lambda_2$ are the eigenvalues of the matrix $M$ in (\ref{4.1}).}

The following corollary is a classification of the expanding
matrix $M$ satisfying (\ref{uno}), $\det M < 0$ and $(\ref{4.1})$.

{\bf Corollary D.} {\em Let $M$ be an expanding  $2 \times
2$-matrix satisfying (\ref{uno}) with $\det M < 0$. Then $M$
satisfies condition $(\ref{4.1})$ if and only if trace $M=0$ and
$\det M = -n$. Especially $$ M^2=nI.$$}

Moreover, they give a classification of the expanding  matrix $M$
satisfying (\ref{uno}),  $\det M > 0$ and $(\ref{4.1})$. They
write the following theorem with the restriction to the case that
the eigenvalues are complex numbers because in other way, by
$(\ref{4.1})$, $M$ is a diagonal matrix. Here, $l$ will denote the
minimal index for which $(\ref{4.1})$ holds, i.e., they neglect
the trivial cases generated by powers of $(\ref{4.1})$.

{\bf Theorem E.} {\em Let $M$ be an expanding  $2 \times 2$-matrix
 satisfying (\ref{uno}) with $\det M > 0$ and eigenvalues
$\lambda_i \notin \mathbb{R}$, $i=1,2$. Then condition
$(\ref{4.1})$ can only hold for $l=3,4,6,8$ and $12$. These cases
can be classified as follows:
\begin{eqnarray*}
M^3=nI & \textrm{if and only if  \ }  \textrm{trace } M = -n^{1/3}
 & \textrm{and }
\det M=n^{2/3},\\
M^4=nI & \textrm{if and only if  \ }  \textrm{trace } M = 0
 & \textrm{and }
\det M=n^{1/2},\\
M^6=nI & \textrm{if and only if  \ }  \textrm{trace } M = n^{1/6}
 & \textrm{and }
\det M=n^{1/3},\\ M^8=nI & \textrm{if and only if  \ }
(\textrm{trace } M)^2 = 2n^{1/4}
 & \textrm{and }
\det M=n^{1/4},\\
 M^{12}=nI & \textrm{if and only if  \ }
(\textrm{trace } M)^2 = 3n^{1/6}
 & \textrm{and }
\det M=n^{1/6}.\\
\end{eqnarray*}}

Moreover,  the following factorization of expanding integer matrices
for some particular cases appears in \cite{DLN:97}.

{\bf Theorem F.} {\em Every  expanding  $2 \times 2$-matrix $M$
 satisfying $(\ref{uno})$ with $\det M =-2^s$, $s \in \{ 1,2,3\}$, which satisfies
  $(\ref{4.1})$ possesses a factorization
  $$
M =ADPA^{-1},
  $$
  where $A$ is a unimodular matrix, $P$ a permutation matrix
  -- that is, $(P\mathbf{x})_i=x_{\pi_i}$ for some permutation $\pi$ of
  $\{1,...,d \}$ and for all $\mathbf{x} \in \mathbb{R}^d$, --
and $D$ is a diagonal matrix with entries $d_i \in \mathbb{Z}$
along the diagonal satisfying
$|d_id_{\pi(i)}...d_{\pi^{d-1}(i)}|>1$ for each $i=1,...,d$.}

Furthermore, in \cite{DGL:97}, the following lemma is proved.

{\bf Lemma G.} {\em Suppose that $M$ is an  expanding $d \times
d$-matrix satisfying $(\ref{uno})$ with the property
 $$
M^d= \pm 2 I.
 $$
If there exists a representative $e \in \mathbb{Z}^d /
M\mathbb{Z}^d$, so that the matrix $(e, Me,...,M^{d-1}e)$ is
unimodular, then $M$ possesses the factorization
  $$
M =AS\Pi A^{-1},
  $$
  where $A \in SL(d,\mathbb{Z})$, $S=diag(\pm2,\pm1,...,\pm1)$
  and $\Pi$ is an irreducible permutation matrix. }


Some results for interpolating scaling functions, although not
that closely related to our topic, can also be found in \cite{DM,
DMT}; in particular, the examples in \cite{DMT} cover MRA even for
some non-selfadjoint matrices, like the quincunx matrix.

\section{Basic notions}\label{basics}

 As a general reference
regarding linear algebra, we refer to \cite{Halmos} and \cite{HK:}.
For further use, and to fix
notation, let us briefly cover some basic facts.\\
Given $r>0$, we denote  $Q_{r}= \{\mathbf{x} \in \RR^{d}~:~ |
x_{i}|< r, \forall i=1,...,d\}$ the cube of side length $2r$ with
the center in the origin.\\
Given a map $A$, we write $d_{A}=|\det A |$. If $A$ is a matrix of
an expansive linear map, then obviously $d_A>1$. The volume of any
measurable set $S$ changes under $A$ according to
$|AS|_d=d_A|S|_d$.\\
A subspace $W \subset \Rd$ is called an \emph{invariant subspace}
under $A$ if $AW \subset W$. As is usual, $W^{\bot}$ is is called
the \emph{orthogonal complement} of $W$ with respect to the
canonical inner product on $\Rd$.
The orthogonal projection of $\mathbf{w}$ onto $W$ is
$P_{W}(\mathbf{w}):=\mathbf{u}$. 

Let $W_1, W_2$ be  vector spaces then $W_1 \oplus W_2$ is the direct
sum of $W_1$ and $W_2$.
 If $A_{\mu} : W_{\mu} \To W_{\mu}$,
$\mu = 1,2$ are linear maps, then we denote by $A_1 \otimes A_2$ the
map on $W_1 \oplus W_2$ defined for any $\mathbf{w_{\mu}} \in
W_{\mu}$, $\mu = 1,2$ as $(A_1 \otimes A_2)(\mathbf{w_1}+
\mathbf{w_2})= A_1
\mathbf{w_1}+ A_2 \mathbf{w_2}$.\\
 If $W$
is an Euclidean space and $A:W \longrightarrow W$ is a linear map,
then $A^*$ will be the \emph{adjoint} of $A$.
$A$ is a self-adjoint map if $A = A^*$.
Let $A_1,A_2:W \longrightarrow W$ be two linear maps. $A_1$ and
$A_2$ are said to be \emph{simultaneously diagonalizable} (see
\cite{HK:}, p. 177) if there exists a basis
$\mathbf{u}_1,...,\mathbf{u}_d$ of $W$ such that $\mathbf{u}_l$, $l
= 1,...,d,$ are eigenvectors of both $A_1$ and $A_2$.\\
The \emph{Spectral Theorem} for self adjoint maps (see \cite[Theorem
1, p.156]{Halmos}) tells us that for any self-adjoint linear map $A$
on $\RR^d$, if $\beta_1<\dots<\beta_k $ are  all the \emph{distinct}
eigenvalues of $A$ with respective multiplicities $m_1,\dots,m_k$,
then for each $i = 1,...,k,$ there exists an orthonormal basis
\[
\mathbf{u}_{m_0+...+m_{i-1}+1},...,
\mathbf{u}_{m_0+...+m_{i-1}+m_{i}},
\]
where $m_0=0$, for the subspace $U_i$ of all  eigenvectors
associated with the eigenvalue $\beta_i$, moreover, then
$\RR^d=U_1\oplus\dots\oplus U_k$ with
\[
U_i= [\mathbf{u}_{m_0+...+m_{i-1}+1},...,
\mathbf{u}_{m_0+...+m_{i-1}+m_{i}}]
\]
being mutually orthogonal, invariant subspaces. Furthermore, we can
then write  $A=\otimes_{i=1}^{k}A_i$ where $A_i:=A|_{U_i}$ are
homothetic transformations $\mathbf{x}\to\beta_i \mathbf{x}$,
$\forall \mathbf{x} \in U_i,$ $i=1,\dots,k$.\\
For a general linear map $M$ on $\RR^d$, one can similarly find a
decomposition $\RR^d =U_1\oplus\dots\oplus U_k$ of invariant
subspaces, which, however, is not necessarily be an
orthogonal decomposition, see \cite[Theorem 2, p.113]{Halmos}.\\
Recall that a linear map $ A:\RR^{d} \longrightarrow \RR^{d}$ is
called \emph{positive} map if it is self-adjoint and all its
(necessarily real) eigenvalues are also positive.

Let $J:\RR^{d}\longrightarrow \RR^{d}$ be a positive map having a
diagonal matrix $J$ and let $\lambda_1,..., \lambda_d \in
[0,\infty)$ be the elements in the diagonal. Then, if
$A:\RR^{d}\longrightarrow \RR^{d}$ is a linear map such that $A= CJ
C^{-1}$ where $C$ is a $d \times d$ invertible matrix, then the
powers $A^t$, $t\in \RR$ are defined as $A^t= CJ^t C^{-1}$, where
$J^t$ is a diagonal matrix with elements $\lambda_1^t,...,
\lambda_d^t$ in the diagonal.

\section{Properties of sets having $0$ as a point of $A$-density}

The next monotonicity property is clear.
\begin{prop} \label{propo106}
Let $A:\RR^{d} \longrightarrow \RR^{d}$ be an expansive linear map.
Let $E, F \subset ~\RR^d$ be measurable sets such that $E \subset F$
and $E \in \mathcal{E}_{A}$. Then $F \in \mathcal{E}_{A}$.
\end{prop}

In the following  propositions we give different equivalent
conditions for the origin to be a  point of $A$-density for a
measurable set $E \subset \RR^d$. Put $E^c:= \RR^d \setminus E$.

\begin{prop}\label{proposicion1} Let $A:\RR^{d}
\longrightarrow \RR^{d}$ be an expansive linear map. Let $E \subset
\RR^{d}$ be a measurable set. Then for any $r>0$ the following four
conditions are equivalent:
\begin{enumerate}
\item [($\mbox{i}$)] \begin{description}
\begin{equation} \lim_{j\longrightarrow \infty}\frac{|
E\bigcap A^{-j}B_{r}|_{d}}{| A^{-j}B_{r}|_d}=1; \label{AD1} 
\end{equation} \end{description}
\item [($\mbox{ii}$)] \begin{description}
\begin{equation} \lim_{j\longrightarrow \infty} |
A^{j}E\bigcap B_{r}|_{d}=| B_{r}|_{d}; \label{AD3} 
\end{equation}\end{description} \item [($\mbox{iii}$)] \begin{description}
\begin{equation} \lim_{j\longrightarrow \infty}\frac{|
E^c\bigcap A^{-j}B_{r}|_{d}}{| A^{-j}B_{r}|_d}=0;
\label{AD2}
\end{equation}\end{description}  \item [($\mbox{iv}$)]
\begin{description}
\begin{equation} \lim_{j\longrightarrow \infty} |
A^{j}E^c\bigcap B_{r}|_{d}=0. \label{AD4} 
\end{equation} \end{description}
\end{enumerate}
\end{prop}
\pf  $(\mbox{i}) \Longleftrightarrow (\mbox{ii}) $ This is a
direct consequence of the fact that for any $r>0$, and for any $j
\in \mathbb{N}$,
\[
\frac{| E\bigcap A^{-j}B_{r}|_{d}}{| A^{-j}B_{r}|_{d}} =\frac{|
A^{j}E\bigcap B_{r}|_{d}}{| B_{r}|_{d}}.
\]
$(\mbox{i}) \Longleftrightarrow (\mbox{iii}) $ Obviously, for any
$r>0$, and for any $j \in \mathbb{N}$,
\[
1= \frac{| E \bigcap A^{-j} B_r |_d}{| A^{-j} B_r |_d} +\frac{| E^c
\bigcap A^{-j} B_r |_d}{| A^{-j} B_r |_d}.
\]
$(\mbox{iii}) \Longleftrightarrow (\mbox{iv})$ This follows since
for any $r>0$, and for any $j \in \mathbf{N}$,
\[
\frac{| E^c\bigcap A^{-j}B_{r}|_{d}}{| A^{-j}B_{r}|_{d}} =\frac{|
A^{j}E^c\bigcap B_{r}|_{d}}{| B_{r}|_{d}}. \qed
\]

\begin{cor}\label{c:trivi} In order to $E\in \EA$, the validity of
any of the above conditions $(i)-(iv)$, but required for all
$r>0$, are necessary and sufficient.
\end{cor}

Two sets are termed \emph{essentially disjoint}, if their
intersection is of measure zero.

\begin{cor}\label{c:disjoint} For any expansive map $A$ and two
sets $E,F\subset \Rd$, which are essentially disjoint, at most one
of the sets can belong to $\EA$.
\end{cor}
\pf Assume, e.g., $E\in\EA$. Note that $F\in\EA$ if and only if
$\widetilde{F}:=F\setminus(E\cap F)\in\EA$, since deleting the
measure zero intersection does not change the measures, hence
neither the limits in the definition of $\EA$. But
$\widetilde{F}\subset E^c$, and $E\in \EA$ entails that the limits
(iii) and (iv) in Proposition \ref{proposicion1} are zero, hence
$E^c\notin \EA$. Obviously (or by the monotonicity formulated in
Proposition \ref{propo106}), then neither $\widetilde{F}\subset E^c$
can belong to $\EA$. Whence $F\notin\EA$. \qed

\begin{prop} \label{proposicion0} Let $A:\RR^{d}\longrightarrow
\RR^{d}$ be an expansive linear map. Let  $E \subset \RR^{d}$ be a
measurable set and assume that for a certain $r_{0}>0$ some (and
hence all) of conditions (i)-(iv) of Proposition \ref{proposicion1}
are satisfied. Then $E \in \mathcal{E}_{A}$. Conversely, if for any
$r_0>0$ any of the conditions (i)-(iv) of Proposition
\ref{proposicion1} fails, then $E \notin \mathcal{E}_{A}$.
\end{prop}
\pf Let $r \in \RR$ and $0 < r < r_0$, and let  $j \in
\mathbf{N}\setminus \{0\}$, then $A^{-j}B_{r} \subset
A^{-j}B_{r_{0}}$, hence by condition $(iii)$
\[
\frac{| E^{c} \bigcap A^{-j}B_{r}|_{d}}{| A^{-j}B_{r}|_{d}} \leq
 (\frac{r_{0}}{r})^{d}\frac{| E^{c} \bigcap
A^{-j}B_{r_{0}}|_{d}}{| A^{-j}B_{r_{0}}|_{d}} \To 0, \textrm{ when
$j \To +\infty$}.
\]

Now let  $r \in \RR$ and $ r>r_{0},$ and let $j \in
\mathbf{N}\setminus \{0\}$. As the map $A$ is  an expansive map,
$\exists m=m(r) \in \mathbf{N}$ such that $B_{r} \subset
A^{m}B_{r_{0}}$. Then similarly to the above
\[
\frac{| E^{c} \bigcap A^{-j}B_{r}|_{d}}{| A^{-j}B_{r}|_{d}} \leq
d_{A}^{m} \frac{| E^{c} \bigcap A^{-j+m}B_{r_{0}}|_{d}}{|
A^{-j+m}B_{r_{0}}|_{d}} \To 0, \textrm{ when $j \To + \infty$.}
\qed
\]

\begin{prop} \label{porpo3} Let $A: \RR^{d} \longrightarrow
\RR^{d}$ be an expansive linear map, and let  $E \subset \RR^{d}$ be
a measurable set. Assume that $K \subset \RR^d$ is another
measurable set, and that there exist $r_1, r_2$ where
$0<r_1<r_2<\infty$ such that $B_{r_1} \subset K \subset B_{r_2}$.
Then $E \in \mathcal{E}_{A}$ if and only if
\begin{equation}
\lim_{j\longrightarrow \infty}\frac{| E\bigcap A^{-j}K |_{d}}{|
A^{-j}K |_{d}} =1,
\end{equation}
or equivalently,
\begin{equation} \lim_{j\longrightarrow \infty}\frac{|
E^c \bigcap A^{-j}K|_{d}}{| A^{-j}K|_{d}} =0. \label{K}
\end{equation}
\end{prop}
\pf $\Longrightarrow)$ In view of the condition $B_{r_1} \subset K
\subset B_{r_2}$ we have
\begin{eqnarray*}
\limsup_{j\longrightarrow \infty}\frac{| E^c\bigcap A^{-j}K|_{d}}{|
A^{-j}K|_{d}}& \leq & \limsup_{j\longrightarrow \infty}\frac{|
E^c\bigcap A^{-j}B_{r_2}|_{d}}{| A^{-j}B_{r_1}|_{d}}\\ & \leq &
(\frac{r_2}{r_1})^d \lim_{j\longrightarrow \infty}\frac{| E^c\bigcap
A^{-j}B_{r_2} |_{d}}{| A^{-j}B_{r_2}|_{d}} =0,
\end{eqnarray*}
because $E \in \mathcal{E}_A.$

$\Longleftarrow)$ Again, by assumption we have
\[
\limsup_{j\longrightarrow \infty}\frac{| E^c\bigcap
A^{-j}B_{r_1}|_{d}}{| A^{-j}B_{r_1}|_{d}} \leq (\frac{r_2}{r_1})^d
\lim_{j\longrightarrow \infty}\frac{| E^c\bigcap A^{-j}K|_{d}}{|
A^{-j}K|_{d}}=0,
\]
using now (\ref{K}). Finally, Proposition \ref{proposicion0} tells
us that the origin is a point of $A$-density for $E$. \qed

\begin{lem} \label{lemma4}
Let $A:\RR^{d}\longrightarrow \RR^{d}$ be an expansive linear map.
Assume that $Y \subset \RR^d$, $Y \cong \RR^p$, $1 \leq p < d$, is
an invariant subspace under $A$, and that also $Y^{\bot}$ is an
invariant subspace under $A$. Let  $E \subset \RR^{d}$ be a
measurable set of the form $E:=Y + F$, where $F \subset Y^\bot$.
Then $E \in \mathcal{E}_{A}$ if and only if $F \in
\mathcal{E}_{A|_{Y^{\bot}}}$.
\end{lem}
\pf As $\RR^d= Y \oplus Y^{\bot}$, and moreover the subspaces $Y$
and $Y^{\bot}$ are invariant subspaces under $A$, we can write $A=
A|_{Y} \otimes A|_{Y^{\bot}}$.

We put $K:=K_1 + K_{2}$, where $K_1:= \{ \mathbf{y} \in Y : \|
\mathbf{y} \| \leq 1 \}$, and $K_{2}:= \{ \mathbf{y} \in Y^{\bot}
: \| \mathbf{y} \| \leq 1 \}$. Observe that $K$ satisfies the
conditions of Proposition \ref{porpo3}.

With this notation, given $j \in \mathbf{N}$ we arrive at
\begin{eqnarray*}
E\bigcap A^{-j}K & = & (Y + F)\bigcap (A|_{Y} \otimes
A|_{Y^{\bot}})^{-j}(K_1 + K_{2}) \\  & = & (Y + F)\bigcap
((A|_{Y})^{-j}K_1 + (A|_{Y^{\bot}})^{-j}K_{2}) \\  & = &
(A|_{Y})^{-j}K_1 + (F \bigcap (A|_{Y^{\bot}})^{-j}K_{2}).
\end{eqnarray*}
As the summands are subsets of $Y$ and $Y^{\bot}$, respectively,
this last sum is also a direct sum. Hence we are led to
$$
\frac{| E\bigcap (A^{-j}K)|_{d}}{| A^{-j}K|_{d}}= \frac{| F
\bigcap (A|_{Y^{\bot}})^{-j}K_{2}|_{d-p}}{|
(A|_{Y^{\bot}})^{-j}K_{2}|_{d-p}}.
$$
Taking limits and applying Proposition \ref{porpo3} we conclude the
proof. \qed

\begin{lem} \label{lemma6}
Let $A_1,A_2:\RR^{d}\longrightarrow \RR^{d}$ be expansive linear
maps and assume that $W \subset \RR^d$ is a subspace of $\RR^d$ such
that both $W$ and $W^{\bot}$ are invariant subspaces under both
$A_1$ and $A_2$. If $\mathcal{E}_{A_1}= \mathcal{E}_{A_2}$ then
$\mathcal{E}_{A_1|_{W}}= \mathcal{E}_{A_2|_{W}}$.
\end{lem}
\pf We consider the cylindrical sets $E = F + W^{\bot}$. According
to Lemma \ref{lemma4} we know that $E \in \mathcal{E}_{A_{\mu}}
\Longleftrightarrow F \in \mathcal{E}_{A_{\mu}|_{W}}$, $\mu =
1,2.$ Therefore, the lemma follows. \qed

\begin{lem} \label{lema1}
Let $A,A':\RR^{d} \longrightarrow \RR^{d}$ be expansive linear maps
and suppose that there is a  linear map $C:\RR^{d} \longrightarrow
\RR^{d}$  with $d_{C}>0$, such that $A'=C^{-1}AC$. Moreover, let $E
\subset \RR^{d}$, $| E |_{d}>0$, be a measurable set. Then $E \in
\mathcal{E}_{A} $ if and only if $C^{-1}E \in \mathcal{E}_{A'}$,
i.e. $\EA=C\EE_{A'}$.
\end{lem}
\pf {$A$, $A'$, $C$ and $C^{-1}$ are invertible linear maps, thus
we have that for any $j \in \mathbf{N} \setminus \{0\}$,
\begin{equation}
\frac{| (C^{-1}E)^{c}\bigcap A'^{-j}B_{1}|_{d}}{|
A'^{-j}B_{1}|_{d}} = \frac{| E^{c}\bigcap A^{-j}CB_{1}|_{d}}{|
A^{-j}CB_{1}|_{d}}. \label{A2}
\end{equation}
Moreover, as $C$ is an invertible linear map, there exists $0 <
r_1 < r_2 < \infty $ such that $B_{r_1} \subset CB_{1} \subset
B_{r_2}$. Therefore, the statement follows from (\ref{A2}) and
Proposition \ref{porpo3}. \qed

A direct consequence of  Lemma \ref{lema1} is the following
\begin{cor}  \label{coro}
Let $A_1,A_2:\RR^{d} \longrightarrow \RR^{d} $ be simultaneously
diagonalizable expansive linear maps. If $| \lambda_i^{(1)} | = |
\lambda_i^{(2)} | $, $i = 1,...,d$, where $\lambda_i^{(\mu)}$,
$\mu=1,2$,  $i = 1,...,d$ are the eigenvalues of $A_{\mu}$,
$\mu=1,2$, then
\[
\mathcal{E}_{A_1}= \mathcal{E}_{A_2}.
\]
\end{cor}
\pf As $A_1$ and $A_2$ are simultaneously diagonalizable, there
exists a linear map $C:\RR^{d} \longrightarrow \RR^{d} $ with
$d_{C}>0$, such that $A_{\mu}=C^{-1}J_{\mu}C$, $\mu= 1,2$. From
Lemma \ref{lema1}, we know that
\[
\mathcal{E}_{A_1}= \mathcal{E}_{A_2} \Longleftrightarrow
\mathcal{E}_{J_1}= \mathcal{E}_{J_2}.
\]
Finally, $\mathcal{E}_{J_1}= \mathcal{E}_{J_2}$ is true because from
$| \lambda_i^{(1)} | = | \lambda_i^{(2)} | $, $i = 1,...,d$ it
follows that for any $j \in \mathbb{Z} $ and for any $r>0$ we have
$J_1^{j}B_{r}=J_2^{j}B_{r}$. \qed

\section{Some particular cases}

\begin{lem} \label{eje3}
Let $A:\RR^{2}\longrightarrow \RR^{2}$ be a diagonal, positive,
expansive linear map with the corresponding matrix
\[
A= \left (
\begin{array}{ccc}
  \lambda_{1} & 0 \\ 0 &  \lambda_{2}
\end{array}  \right ), \qquad  \lambda_{1},
\lambda_{2} \in \RR, ~ 1 <  \lambda_{1}, \lambda_{2}.
\]
Let for any $\alpha>0$ $E_{\alpha} \subset \RR^{2}$ be the set
\[
E_{\alpha}= \{(x_{1},x_{2}) \in \RR^{2}: | x_{2} | \geq | x_{1 } |
^{\alpha} \}.
\]
Denote $\alpha_{1,2}:=\alpha_{1,2}(\lambda_1,\lambda_2):= \log
\lambda_2/\log\lambda_1$. Then $E_{\alpha} \in \mathcal{E}_{A}$ if
and only if $\alpha > \alpha_{1,2}$.
\end{lem}
\pf For any $j \in \mathbb{N} \setminus \{0\}$, and because of the
symmetry of the sets $E_{\alpha}^c$ and $ A^{-j}Q_{1}$
\begin{eqnarray} \label{alphacalc}
| E_{\alpha}^{c} \bigcap A^{-j}Q_{1}|_{2} & = 4\int_{0}^{
\lambda_{1}^{-j}} \int_0^{\lambda_2^{-j}}
{\mathbf{1}}_{\{x_2<x_1^{\alpha}\}} d x_2 dx_{1}= 4\int_{0}^{
\lambda_{1}^{-j}} \min(x_{1}^{\alpha},\lambda_2^{-j}) dx_{1},
\end{eqnarray}
for any value of $\alpha>0$. Let us consider first the boundary
case $\alpha=\alpha_{1,2}$. Then $x_1^{\alpha_{1,2}}\leq
\lambda_1^{-j\alpha_{1,2}}=\lambda_2^{-j}$, hence the minimum is
just $x_1^{\alpha_{1,2}}$, and we get
\[
| E_{\alpha_{1,2}}^{c} \bigcap A^{-j}Q_{1}|_{2} = 4\int_{0}^{
\lambda_{1}^{-j}} x_{1}^{\alpha_{1,2}} dx_{1}= 4\frac{
\lambda_{1}^{-j (\alpha_{1,2} + 1)}}{\alpha_{1,2} + 1}.
\]
Therefore,
\[
\frac{| E_{\alpha_{1,2}}^{c} \bigcap A^{-j}Q_{1}|_{2}}{|
A^{-j}Q_{1}|_{2}} = \frac{ \lambda_{1}^{j} \lambda_{2}^{j}}{
\lambda_{1}^{j (\alpha_{1,2} + 1)}(\alpha_{1,2} + 1)} =
\frac{1}{\alpha_{1,2}+1} \left(
\frac{\lambda_{2}}{\lambda_{1}^{\alpha_{1,2}}} \right)^j
=\frac{1}{\alpha_{1,2}+1}
\]
in view of $\lambda_2=\lambda_1^{\alpha_{1,2}}$. The quotient of the
measures on the left being constant, obviously the limit is positive
but less than $1$, hence by Proposition \ref{porpo3} and Proposition
\ref{proposicion1} (i) and (iii) neither $E_{\alpha_{1,2}}$, nor its
complement $E^c_{\alpha_{1,2}}$ can belong to $\EA$. \\
Note that when $\alpha > \alpha_{1,2}$, then $x_1^{\alpha}\leq
x_1^{\alpha_{1,2}}$ (as $x_1<1$), hence in (\ref{alphacalc})
 the
minimum is again $x_1^{\alpha}$. Therefore, a very similar
calculation as above yields
\[
\lim_{j \To \infty } \frac{| E_{\alpha}^{c} \bigcap
A^{-j}Q_{1}|_{2}}{| A^{-j}Q_{1}|_{2}} = \frac{1}{\alpha+1}
\lim_{j\To\infty} \left( \frac{\lambda_{2}}{\lambda_{1}^{\alpha}}
\right)^j = 0,
\]
because now we have $ \lambda_{2}/\lambda_{1}^{\alpha}=
\lambda_{1}^{\alpha_{1,2}-\alpha} < 1$. Whence Proposition
\ref{porpo3} and Proposition \ref{proposicion1} (iii) now gives
$E_\alpha\in \EA$. \\
Finally, let $\alpha<\alpha_{1,2}$. Observe that the coordinate
changing isometry of $\RR^2$ provides a symmetry for our subject:
changing the role of the coordinates we can consider now
$\widetilde{E}_{\beta}:= \{ (x_1,x_2)\in\RR^2~:~ |x_1|\geq
|x_2|^{\beta}\}$. Then obviously $E^c_{\alpha}=
\textrm{int}~\widetilde{E}_{1/\alpha}
\subset\widetilde{E}_{1/\alpha}$, and
$\widetilde{\alpha}_{1,2}=\alpha_{2,1}=\log\lambda_1/\log\lambda_2
=1/\alpha_{1,2}$, hence from the previous case and Proposition
\ref{propo106} we obtain $E^c_{\alpha} \in \EA$. But then
$E_{\alpha} \notin \EA$. That finishes the proof of the Lemma.
\qed

\begin{lem}\label{eje4}
Let $A:\RR^{d}\longrightarrow \RR^{d}$ be a positive expansive
linear map. With the notation in \S \ref{basics}, given $\delta>0$,
we define the measurable set
\begin{eqnarray*}\label{Gdelta}
G_{\delta}:&=& \{ \mathbf{x}=\mathbf{y}+\mathbf{z}~:~
\mathbf{y}\in U_1,~ \mathbf{z} \in U_1^\bot, \|\mathbf{z}\|<
\delta \|\mathbf{y}\|\}  \\ & = & \{
\mathbf{x}=\sum_{i=1}^k\mathbf{y}_i ~:~ \mathbf{y}_i\in U_i,
i=1,\dots,k, ~\| \mathbf{y_2}+...+ \mathbf{y_k}\|< \delta \|
\mathbf{y_1}\| \}.
\end{eqnarray*}
Then in case $\dim U_1 <d$, i.e. when not all the eigenvalues are
equal to $\beta_1$, we have $G_{\delta}\in \mathcal{E}_{A}$.
\end{lem}
\pf Clearly, $|B_1\bigcap G_{\delta}|_d = \int_{B_1}
\mathbf{1}_{G_\delta} d\mathbf{x}$, and $\mathbf{1}_{G_\delta} \To
\mathbf{1}_{\Rd}$ a.e. when $\delta\to~\infty$, so by the Lebesgue
dominated convergence theorem we conclude
\begin{equation}\label{Lebesgue}
\lim_{\delta \To \infty} | B_1\bigcap G_{\delta} |_d =  | B_1 |_d.
\end{equation}
Next we prove that for any given $\delta>0$,
$G_{\frac{\beta_2}{\beta_1}\delta}\subset A G_{\delta}$. We can
write $A G_{\delta}$ as
\begin{eqnarray*}
A G_{\delta}&= &\{ \sum_{i=1}^k \beta_i\mathbf{y}_i ~:~ \mathbf{y}_i
\in U_i, i=1,...,k, \sum_{i=2}^k \| \mathbf{y}_i\parallel^2 <
\delta^2 \|\mathbf{y}_1 \|^2 \}
\\ &=&\{\sum_{i=1}^k \mathbf{z}_i ~:~ \mathbf{z}_i\in U_i,
i=1,...,k, \sum_{i=2}^k  \frac{1}{\beta_i^2}  \| \mathbf{z}_i\|^2 <
\delta^2 \frac{1}{\beta_1^2}\| \mathbf{z}_1\|^2 \}.
\end{eqnarray*}
Let $\mathbf{x}\in G_{\frac{\beta_2}{\beta_1}\delta}$. Then,
\[ \mathbf{x}=\mathbf{z}_1+...+
\mathbf{z}_k \textrm{ such that } \frac{1}{\beta_2^2}\|
\mathbf{z}_2\|^2+...+\frac{1}{\beta_2^2} \|\mathbf{z}_k\|^2 <
\frac{1}{\beta_1^2}\delta^2 \|\mathbf{z}_1\|^2,
\]
and as $\beta_i > \beta_2$, $i=3,...,k$, then
\[
\frac{1}{\beta_2^2}\| \mathbf{z}_2\|^2+...+\frac{1}{\beta_k^2}
\|\mathbf{z}_k\|^2 \leq \frac{1}{\beta_2^2}\|
\mathbf{z}_2\|^2+...+\frac{1}{\beta_2^2} \|\mathbf{z}_k\|^2 <
\frac{1}{\beta_1^2}\delta^2 \|\mathbf{z}_1\|^2.
\]
Hence we arrive at $\mathbf{x}\in A G_{\delta}$ proving
$G_{\frac{\beta_2}{\beta_1}\delta}\subset A G_{\delta}$ indeed. If
now we iterate this and use (\ref{Lebesgue}), we infer
\[
\lim_{j \To \infty} | B_1\bigcap A^jG_{\delta} |_d \geq \lim_{j
\To \infty} | B_1\bigcap G_{(\frac{\beta_2}{\beta_1})^j\delta}|_d
=  | B_1 |_d,
\]
so by Proposition \ref{proposicion1} (ii) and Proposition
\ref{proposicion0} we get $G_{\delta}\in \mathcal{E}_{A}$. \qed

\begin{lem}  \label{eje5}
Let $A:\RR^{d}\longrightarrow \RR^{d}$ be a positive expansive
linear map, and similarly to \S \ref{basics} let the different
eigenvalues be listed as $1< \beta_1<...< \beta_k$, $U_1 \subset
\RR^d$ being the eigenspace belonging to $\beta_1$. Moreover, let $V
\subset U_1^{\bot}$ be any subspace of $\Rd$ orthogonal to $U_1$,
and write $W:=(U_1\oplus V)^{\bot}$. We finally set for any
$\delta>0$
\begin{eqnarray*}
F_{\delta} = \{ \mathbf{x} = \mathbf{u}+\mathbf{v}+\mathbf{w}~:~
\mathbf{u}\in U_1, ~\mathbf{v}\in V, ~\mathbf{w}\in (U_1 \oplus
V)^{\bot},~ ~ \|\mathbf{v} \|< \delta \|\mathbf{u} \|\}.
\end{eqnarray*}
Then  $F_{\delta} \in \mathcal{E}_{A}$.
\end{lem}
\pf We can combine Proposition \ref{propo106} and  Lemma
\ref{eje4}, because $G_{\delta}$
 is contained in
$F_{\delta}$. \qed

\section{The Main Result}

\begin{thm}  \label{thm}
Let $A_1,A_2:\RR^{d}\longrightarrow \RR^{d}$ be positive expansive
linear maps. Then $\mathcal{E}_{A_1}=\mathcal{E}_{A_2}$ if and only
if $\exists t>0$ such that
\[
A_1^t=A_2.
\]
\end{thm}

For the proof of Theorem \ref{thm}, we first settle the case of
diagonal matrices in the following lemma. After that, we will
apply the spectral theorem to prove even the general case.

\begin{lem} \label{log1=log2}
Let $A_1,A_2:\RR^{d}\longrightarrow \RR^{d}$ be positive diagonal
expansive linear maps with the corresponding matrices
\[ A_{\mu}= \left (
\begin{array}{ccccccccc}
  \lambda_{1}^{(\mu)} & 0  & 0 & ... & 0  \\
  0 &  \lambda_{2}^{(\mu)} & 0 & ... & 0  \\
  .. & .. & .. & ... & .. \\
  0 & 0 & 0 & ... & \lambda_{d}^{(\mu)}
\end{array}  \right ), \] where   $
\lambda_{i}^{(\mu)} \in \RR$,  $1 <  \lambda_{1}^{(\mu)} \leq
\lambda_{2}^{(\mu)} \leq ... \leq  \lambda_{d}^{(\mu)}$, for $\mu =
1,2$. Then $\mathcal{E}_{A_1} = \mathcal{E}_{A_2} $ if and only if
$\exists t>0$ such that
\[
(A_1)^t=A_2.
\]
\end{lem}
\pf $\Rightarrow)$ For an indirect proof, we assume that it is false
that $\exists t>0$ such that$ (A_1)^t=A_2 $. Then $\exists i,l \in
\{1,...,d\},$ $i<l$, such that
$(\lambda_{i}^{(1)})^{t_1}=\lambda_{i}^{(2)}$ and
$(\lambda_{l}^{(1)})^{t_2}=\lambda_{l}^{(2)}$ with $0< t_1, t_2$ but
$t_1 \neq t_2$, i.e.
\[
t_1 = \frac{\ln  \lambda_{i}^{(2)}  }{ \ln  \lambda_{i}^{(1)}} \neq
\frac{\ln  \lambda_{l}^{(2)}  }{ \ln  \lambda_{l}^{(1)}  }= t_2 ,
\]
or equivalently
\[
\alpha_1:= \frac{\ln  \lambda_{l}^{(1)}  }{ \ln
 \lambda_{i}^{(1)}  } \ne \alpha_2:= \frac{\ln  \lambda_{l}^{(2)}  }{
 \ln \lambda_{i}^{(2)}
} .
\]
Without loss of generality, we can assume that $(1\leq)\alpha_1 <
\alpha_2$. Let $\alpha>0$ and let us define
\[
F := \{ (x_i,x_l)\in \RR^{2}: | x_l | \geq | x_i |^{\alpha} \}
\]
and
$$
E := \{\mathbf{x}= (x_1,...,x_d)\in \RR^{d}~:~ | x_l | \geq | x_i
|^{\alpha},~x_j\in\RR ~(j\ne i,l)~ \}\cong F\oplus\RR^{d-2}.
$$
Then Lemma \ref{lemma4} tells us that $E \in \mathcal{E}_{A_{\mu}} \
\ \mu = 1,2 \Longleftrightarrow   F \in \mathcal{E}_{M_{\mu}} \ \
\mu = 1,2$, where $M_1, M_2:\RR^{2}\longrightarrow \RR^{2}$ are
expansive linear maps with matrices
\[
M_{\mu}= \left (
\begin{array}{ccc}
  \lambda_{i}^{(\mu)} & 0 \\ 0 &  \lambda_{l}^{(\mu)}
\end{array}  \right ),  \ \ \ \ \ \ \mu = 1,2.
\]
However, making use of $\alpha_1 < \alpha_2$, we can choose a value
$\alpha_1 < \alpha < \alpha_2$, and then Lemma \ref{eje3} gives $F
\in \mathcal{E}_{M_{1}} $ but $F \notin \mathcal{E}_{M_{2}}$,
contradicting to the assumption
$\mathcal{E}_{A_{1}}=\mathcal{E}_{A_{2}}$.

$\Leftarrow)$ As $A_2=A_1^t$ if and only if $A_1=A_2^{1/t}$, it
suffices to see that $\mathcal{E}_{A_1} \subset \mathcal{E}_{A_2}$.
So let $E\in\EE_{A_1}$.

Since $A_1$ is a positive, expansive diagonal mapping, obviously
for any $0\leq s <1$ we have $B_1\subset A_1^{s} B_1 \subset A_1
B_1$. Now write, for any $j \in \NN \setminus \{0\}$, the exponent
$tj$ as $tj=l_j-s_j$ with $l_j:=\lceil tj \rceil$,  the least
integer $\geq tj$, and $s_j:=\lceil tj \rceil - tj \in[0,1)$. So
we have
\[
\frac{| E^{c} \bigcap A_2^{-j}B_{1}|_{d}}{| A_2^{-j}B_{1}|_{d}} =
\frac{| E^{c} \bigcap A_1^{-l_{j}+s_j}B_{1}|_{d}}{|
A_1^{-l_{j}+s_j}B_{1}|_{d}} 
\leq d_{A_1}\frac{| E^{c} \bigcap A_1^{-l_{j}+1}B_{1}|_{d}}{|
A_1^{-l_{j}+1}B_{1}|_{d}}.
\]
 Since $\{-l_j+1\}_{j
\in \NN}$ is an integer sequence and $-l_j+1 \to -\infty$ when $j
\to \infty$,  by condition $E\in \EE_{A_1}$, Proposition
\ref{proposicion1} (iii) entails that the right hand side
converges to $0$ with $j\to\infty$, whence
\[
\lim_{j \To \infty}\frac{| E^{c} \bigcap A_2^{-j}B_{1}|_{d}}{|
A_2^{-j}B_{1}|_{d}}=0.
\]
According to  Proposition \ref{proposicion0} this means
$E\in\EE_{A_2}$. \qed

\begin{lem}  \label{lemmaU}
Let $A_1,A_2:\RR^{d}\longrightarrow \RR^{d}$ be positive expansive
linear maps such that $\mathcal{E}_{A_1}=\mathcal{E}_{A_2}$. Then
$\dim U_1^{(1)}\bigcap U_1^{(2)} \geq 1$.
\end{lem}
\pf Assume the contrary, i.e. $U_1^{(1)}\bigcap U_1^{(2)}=\{
\mathbf{0} \}$, hence $V:=U_1^{(1)}+U_1^{(2)}= U_1^{(1)}\oplus
U_1^{(2)}$. Recall that by definition both $U_1^{(1)}$ and
$U_1^{(2)}$ are of dimension at least one, and now $\dim
U_1^{(1)}+\dim U_1^{(2)} = \dim V :=p \leq d$, hence now neither
of them can have full dimension. Without loss of generality we can
assume $V=\RR^p$. First we work in $V$. Denote
$V_{\mu}:=(U_1^{(\mu)})^{\bot}$ (the orthogonal complement
understood within $V$) for $\mu=1,2$.

If $S$ is the unit sphere of $V$,  $S:=\{\mathbf{x}\in V~:~
\|\mathbf{x}\| =1\}$, then by he indirect assumption also the
traces $T_{\mu}:=S\cap U_1^{(\mu)}$ are disjoint for $\mu=1,2$. As
these sets are compact, too, there is a positive distance
$0<\rho:=
\textrm{dist}(T_1,T_2)\leq \sqrt{2}$ between them.\\
Let us fix some parameter $0<\kappa<1$, to be chosen later. Next
we define the sets
\[
K_{\mu}:=\{\mathbf{u}+\mathbf{v}~:~ \mathbf{u}\in U_1^{(\mu)},
\mathbf{v}\in V_{\mu}, \|\mathbf{v}\|\leq \kappa \|\mathbf{u}\| \}
\qquad (\mu=1,2).
 \]
  We claim that these sets are essentially
disjoint, more precisely $K_1\cap K_2=\sn$, if $\kappa$ is chosen
appropriately. So let now $\mu=1$ or $\mu=2$ be fixed, and
consider any $\xx\in K_{\mu}$ with $\gamma:=\|\xx\|\ne 0$, i.e.
$\xx\in K_\mu\setminus\sn$. From the representation of $\xx$ as
the sum of the orthogonal vectors $\uu$ and $\vv$, we get $\|\uu\|
\leq \| \xx \| = \sqrt{\|\uu\|^2+ \|\vv\|^2} \leq \sqrt{\|\uu\|^2
+ \kappa^2 \|\uu\|^2} = \sqrt{1+\kappa^2} \|\uu\|$. We put
$\beta:=\|\uu\|$. Let now $\yy:=(1/\gamma) \xx\in S$ be the
homothetic projection of $\xx$ on $S$. Then
\begin{eqnarray*}
\textrm{dist}(\yy,T_\mu) & \leq & \left\| \yy-\frac{1}{\beta} \uu
\right\| \leq \left\| \yy-\frac{1}{\beta} \xx \right\| + \left\|
\frac{1}{\beta} \xx - \frac{1}{\beta} \uu \right\| \\ & \leq &
\left| 1- \frac{\gamma}{\beta} \right| + \frac1\beta\|\vv\| \leq
(\sqrt{1+\kappa^2}-1) + \kappa < 2\kappa.
\end{eqnarray*}
Therefore, if we choose $\kappa<\rho/4$, then $\yy$ falls in the
$\rho/2$ neighborhood of $T_{\mu}$, whence the homothetic
projections $\yy_\mu$ of elements $\xx_\mu\in K_\mu$, $\mu=1,2$,
can never coincide. But $K_{\mu}$ are cones, invariant under
homothetic dilations, therefore this also implies that $K_1\cap
K_2 \subset \sn$, as we needed.

Let us write $W:=\left( U_1^{(1)}\oplus U_1^{(2)}\right)^{\bot}$.
Now we consider the sets
\begin{eqnarray*}
H_{\mu}: & = & K_{\mu}\oplus W =\{\mathbf{x}+ \mathbf{w}~:~
\mathbf{x}\in K_{\mu}, \mathbf{w}\in W \} \\ & = &
\{\uu+\vv+\mathbf{w}~:~ \uu\in U_1^{(\mu)}, \vv\in V_{\mu},
\mathbf{w}\in W, \|\vv\|\leq \kappa \|\uu\| \} \quad (\mu=1,2),
\end{eqnarray*}
which are also essentially disjoint, as $H_1\cap H_2 =W$ and
$|W|_d=0$ because $\dim W <d$. These sets are exactly of the form
$F_\delta$ in  Lemma \ref{eje5}, thus $H_{\mu} \in ~
\mathcal{E}_{A_{\mu}}$ for $\mu=1,2$. It remains to recall
Corollary \ref{c:disjoint}, saying that essentially disjoint sets
can not simultaneously be elements of the same $\EE_{A_\mu}$, that
is, $H_1\in \EE_{A_1}$ but then $H_2\notin \EE_{A_1}$, and $H_2\in
\EE_{A_2}$, but $H_2\notin \EE_{A_1}$. Here we arrived at a
contradiction with $\mathcal{E}_{A_1}=\mathcal{E}_{A_2}$, which
concludes our proof. \qed

\begin{lem}  \label{lemmaSimul}
Let $A_1,A_2:\RR^{d}\longrightarrow \RR^{d}$ be positive expansive
linear maps such that $\mathcal{E}_{A_1}=\mathcal{E}_{A_2}$. Then
$A_1$ and $A_2$ are simultaneously diagonalizable maps.
\end{lem}
\pf We prove the lemma by induction with respect to the dimension.
Obviously, the lowest dimensional case of $d=1$ is true. Now let $d
\geq 1$ and assume that for any two positive expansive linear maps
$M_1,M_2:\RR^{d}\longrightarrow \RR^{d}$ such that
$\mathcal{E}_{M_1}=\mathcal{E}_{M_2}$, $M_1$ and $M_2$ are
simultaneously diagonalizable. We will prove that the statement is
true for dimension $d+1$. Let $A_1,A_2:\RR^{d+1}\longrightarrow
\RR^{d+1}$ be positive expansive linear maps such that
$\mathcal{E}_{A_1}=\mathcal{E}_{A_2}$. From  Lemma \ref{lemmaU} we
know that there exists a one dimensional subspace, say $[u]$, so
that $[u] \subset U_1^{(1)}\bigcap U_1^{(2)}$.

As $u$ is an eigenvector of the positive self-adjoint linear maps
$A_1$ and $A_2$, $[u]$ is an invariant subspace of both $A_1$ and
$A_2$, and we have that also $[u]^{\bot}$ is an invariant subspace
under both $A_1$ and $A_2$. Hence from Lemma \ref{lemma6}, we
obtain that $\mathcal{E}_{M_1}=\mathcal{E}_{M_2}$ where  $M_{\mu}:
= A_{\mu}|_{[u]^{\bot}}$, $\mu = 1,2.$ Then by hypothesis of
induction we know that the positive expansive linear maps
$M_1,M_2:[u]^{\bot}\longrightarrow [u]^{\bot}$ are simultaneously
diagonalizable maps. Furthermore, as we can write $A_{\mu}=
A_{\mu} |_{[u]} \otimes M_{\mu}$, $\mu = 1,2,$ and $u \in [u]$ is
an eigenvector of $A_1$ and $A_2$, we can conclude that $A_1$ and
$A_2$ are simultaneously diagonalizable maps. \qed

\begin{pf*}{Proof of Theorem \ref{thm}} $\Longleftarrow)$ From the {\em
spectral theorem} we know that there exists a linear map
$C:\RR^{d}\longrightarrow \RR^{d}$ with $d_{C}>0$, such that
$A_1=CJ_1 C^{-1}$ where $J_1:\RR^{d}\longrightarrow \RR^{d}$ is an
expansive linear map with corresponding matrix
\[
J_{1}= \left (
\begin{array}{ccccccccc}
  \lambda_{1}^{(1)} & 0  & 0 & ... & 0  \\
  0 &  \lambda_{2}^{(1)} & 0 & ... & 0  \\
  .. & .. & .. & ... & .. \\
  0 & 0 & 0 & ... & \lambda_{d}^{(1)}
\end{array}  \right ), \ \ \ \ \ \ \  \textrm{  $
\lambda_{i}^{(1)} \in \RR$,  $1 <  \lambda_{1}^{(1)} \leq
\lambda_{2}^{(1)} \leq ... \leq  \lambda_{d}^{(1)}$.}
\]
According to the condition $A_2=A_1^t$ with $t>0$, we can write
the corresponding matrix of the map $A_2$ as $A_2=C(J_1)^t
C^{-1}$.

Lemma \ref{log1=log2} tells us that
$\mathcal{E}_{J_1}=\mathcal{E}_{(J_1)^t}$. Also we have the
equivalence
\[
\mathcal{E}_{J_1}=\mathcal{E}_{(J_1)^t} \Longleftrightarrow
C\mathcal{E}_{J_1}=C\mathcal{E}_{(J_1)^t}.
\]
Finally, Lemma \ref{lema1} implies $C\EE_{J_1}=\EE_{A_1}$ and
$C\EE_{(J_1)^t}=\EE_{A_2}$, hence
\[
C\mathcal{E}_{J_1}=C\mathcal{E}_{(J_1)^t} \Longleftrightarrow
\mathcal{E}_{A_1}=\mathcal{E}_{A_2}.
\]
This concludes the proof of the $\Longleftarrow)$ direction. \\
$\Longrightarrow)$ According to Lemma \ref{lemmaSimul}, there exists
an orthonormal basis for $\RR^d$, $\mathbf{u}_{1}, ...,
\mathbf{u}_{d}$, such that $A_1$ and $A_2$ have a common diagonal
representation matrix $C$ in this basis. More precisely, the linear
map $C:\RR^{d}\longrightarrow \RR^{d}$, has matrix
$C:=(\mathbf{u}_1,\mathbf{u}_2,..., \mathbf{u}_d)$ (where
$\mathbf{u}_l$, $l = 1,...,d$ are column vectors, formed from the
common eigenvectors of $A_1$ and $A_2$), and we can write
$A_{\mu}=CJ_{\mu} C^{-1}$, $\mu = 1,2$, where
$J_{\mu}:\mathbb{R}^{d}\longrightarrow \mathbb{R}^{d}$ are expansive
diagonal linear maps with the corresponding matrices being
\[
J_{\mu}= \left (
\begin{array}{ccccccccc}
  \lambda_{1}^{(\mu)} & 0  & 0 & ... & 0  \\
  0 &  \lambda_{2}^{(\mu)} & 0 & ... & 0  \\
  .. & .. & .. & ... & .. \\
  0 & 0 & 0 & ... & \lambda_{d}^{(\mu)}
\end{array}  \right ), \quad  \textrm{  $
\lambda_{i}^{(\mu)} \in \RR$,  $1 <  \lambda_{1}^{(\mu)} \leq
\lambda_{2}^{(\mu)} \leq ... \leq  \lambda_{d}^{(\mu)}$.}
\]
Note that $d_C=1>0$ because the orthogonality of the column vectors
$\mathbf{u}_l$ ($l=1,\dots,d$).\\
>From Lemma \ref{lema1} we get
\[
\mathcal{E}_{A_1}=\mathcal{E}_{A_2} \Longleftrightarrow
C\mathcal{E}_{J_1}=C\mathcal{E}_{J_2} \Longleftrightarrow
\mathcal{E}_{J_1}=\mathcal{E}_{J_2}.
\]
And finally, Lemma \ref{log1=log2} tells us that
\[
\mathcal{E}_{J_1}=\mathcal{E}_{J_2} \ \ \Longleftrightarrow \ \
\exists t>0 \ \ \textrm{such that} \ \ (J_1)^t=J_2.
\]
Therefore, we can write $A_2= C (J_1)^t  C^{-1}= (A_1)^t$, which
concludes the proof of Theorem \ref{thm}. \qed
\end{pf*}
For a slightly more general result, let now
$A_1,A_2:\mathbb{R}^{d}\longrightarrow \mathbb{R}^{d}$ be
self-adjoint expansive linear maps, without assuming that they are
positive. We now consider the diagonal matrices
\begin{equation}\label{Jmu}
J_{\mu}= \left (
\begin{array}{ccccccccc}
  \lambda_{1}^{(\mu)} & 0  & 0 & ... & 0  \\
  0 &  \lambda_{2}^{(\mu)} & 0 & ... & 0  \\
  .. & .. & .. & ... & .. \\
  0 & 0 & 0 & ... & \lambda_{d}^{(\mu)}
\end{array}  \right ), \ \ \ \ \ \ \   \textrm{  $
\lambda_{i}^{(\mu)} \in \RR$. }
\end{equation}
where $A_{\mu}= C_{\mu} J_{\mu} C_{\mu}^{-1}$, $\mu = 1,2,$ with
some invertible mappings $C_1,C_2:\RR^{d}\longrightarrow \RR^{d}$.
Let us denote
\begin{equation}\label{Jmuprime}
J'_{\mu}= \left (
\begin{array}{ccccccccc}
  | \lambda_{1}^{(\mu)} | & 0  & 0 & ... & 0  \\
  0 &  | \lambda_{2}^{(\mu)} | & 0 & ... & 0  \\
  .. & .. & .. & ... & .. \\
  0 & 0 & 0 & ... & | \lambda_{d}^{(\mu)} |
\end{array}  \right ).
\end{equation}
Then as a consequence of Theorem \ref{thm} and Corollary
\ref{coro}, we can say the following.
\begin{cor}\label{cormain}
Let $A_1,A_2:\mathbb{R}^{d}\longrightarrow \mathbb{R}^{d}$ be
self-adjoint expansive linear maps. Let
$C_1,C_2:\RR^{d}\longrightarrow \RR^{d}$ be such that $A_{\mu}=
C_{\mu} J_{\mu} C_{\mu}^{-1}$, $\mu = 1,2,$ where $J_{\mu}$ are the
respective diagonal maps as in (\ref{Jmu}). Then
$\mathcal{E}_{A_1}=\mathcal{E}_{A_2}$ if and only if $\exists t>0$
such that
\[
A_1^{'t}=A'_2,
\]
where  $A'_{\mu}= C_{\mu} J'_{\mu} C_{\mu}^{-1}$, $\mu = 1,2,$
with $J'_{\mu}$ in (\ref{Jmuprime}).
\end{cor}

\section{Application  to multiresolution analysis}

In this section we study equivalence among expansive matrices
satisfying (\ref{uno}). In general, the problem is still open. We
look for some description of self-adjoint expansive linear maps
$A_1, A_2 :\RR^d \To \RR^d$ satisfying (\ref{uno}), such that
$\mathcal{E}_{A_1}= \mathcal{E}_{A_2}$. Hence we can get
equivalent self-adjoint expansive linear maps for MRA.

Above we obtained that if $A_1, A_2 :\RR^d \To \RR^d$ are
expansive positive linear maps, then $\mathcal{E}_{A_1}=
\mathcal{E}_{A_2}$ if and only if there exists $t>0$ such that
$A_2 = A_1^t$. The general case of self-adjoint maps reduces to
this case according to Corollary~\ref{cormain}, so in the
following discussion we restrict to this case of positive
equivalent mappings. To meaningfully interpret the general
requirement, one assumes $\mathcal{E}_{A_1}= \mathcal{E}_{A_2}$,
-- so according to Theorem \ref{thm} we have $A_2 = A_1^t$, $t>0$
-- and now we look for further properties to ensure (\ref{uno}),
too. So in the following let us assume that (\ref{uno}) is
satisfied by $A_1$ and by $A_2$.
\\
To fix notations we have already settled with choosing $\ZZ^d$ to be
the fundamental lattice for our MRA. Therefore, we can assume that
$A_1$ is written in diagonal form in the canonical basis of $\ZZ^d$
(otherwise considerations should change to the fundamental lattice
spanned by the orthogonal basis of eigenvectors for $A_1$). As a
consequence of $A_2=A_1^t$, also $A_2$ is in diagonal form with
respect to the canonical basis. Therefore, (\ref{uno}) means that we
require these diagonal entries -- eigenvalues of $A_{\mu}$ -- belong
to $\ZZ$, or, actually, to $\NN$ as they are positive matrices.
\\
In case all eigenvalues of $A_1$ are equal, i.e. $\beta_1^{(1)}$,
by $A_2=A_1^t$ we have the same property also for $A_2$, and the
equation we must solve is that
$\beta_1^{(2)}=(\beta_1^{(1)})^t\in\NN$ and $\beta_1^{(1)}\in\NN$
simultaneously. Clearly, with $t:=\log \beta_1^{(2)} / \log
\beta_1^{(1)}$ this can always be solved, so any two integer
dilation matrices define equivalent MRA. Let us remark that in the
thesis \cite{Thesis} there is a complete analysis of equivalence
(with respect to the notion of points of $A$-density) to the
dyadic dilation matrix, among all expansive linear mappings,
self-adjoint or not. However, our focus here is different, as here
we consider, under assumptions of self-adjointness, equivalence of
arbitrary, not necessarily dilation mappings.\\
In the general case when $A_1$ (and hence also $A_2$) are not
dilations, there must be two different entries (eigenvalues) in
the diagonal of $A_1$ (and of $A_2$). As equivalence is hereditary
in the sense that the restricted mappings on eigensubspaces of
$A_{\mu}$ must also be equivalent, we first restrict to the case
of dimension 2.
\\
In dimension 2, we thus assume that $A_1$ has diagonal elements
$a\ne b$ belonging to $\NN \setminus\{ 0, 1 \}$ and zeroes off the
diagonal, and we would like to know when do we have with some $t>0$
that $a, b, a^t, b^t\in \mathbb{Z}$ (or $\in\NN$). Obviously, if $ t
\in \NN \setminus \{ 0 \}$ then this condition holds for any $a,b\in
\NN$. Also, in case $a$ and $b$ are full $q^{\rm th}$ powers, we can
as well take $t=p/q\in\QQ$ with arbitrary $p\in\NN$. That system of
solutions -- $a=\alpha^q$, $b=\beta^q$, $t=p/q$ with $\alpha, \beta,
p \in \NN$ -- form one
trivial set of solutions for our equivalence.\\
Another trivial set of solutions arises when $b=a^k$ with some $k\in
\NN$. Then it suffices to have $a^t\in \NN$, which automatically
implies $b^t\in\NN$. More generally, if $b=a^{k/m}$ is a rational
relation between $a$ and $b$, then by the unique prime factorization
we conclude that $a$ is a full $m^{\rm th}$ power and that $b$ is
full $k^{\rm th}$ power, and again we find a system of solutions for
all $t\in\QQ$ of the form $t=\ell/k$.\\
All these trivial solutions can be summarized as cases of rational
relations between $a$, $b$ and $t$: once there is such a relation,
one easily checks, if the respective matrix entries really become
integers. So we find that systems of trivial solutions do exist if
either $t$ is rational, or if $\log a$ and $\log b$ are rationally
dependent (are of rational multiples of each other). We can thus
call these cases the \emph{trivial equivalence} of self-adjoint
expansive linear maps with respect to MRA construction. These
explain the next definition.

\begin{defn}\label{def:trivieq} Let $A_\mu$ ($\mu=1,2$) be two
self-adjoint expansive linear maps, with $A_{\mu}= C_{\mu} J_{\mu}
C_{\mu}^{-1}$, where $J_{\mu}$ are the respective diagonal maps as
in (\ref{Jmu}), and $A'_{\mu}= C_{\mu} J'_{\mu} C_{\mu}^{-1}$, $\mu
= 1,2,$ with $J'_{\mu}$ in (\ref{Jmuprime}) for $\mu = 1,2,$. We say
that $A_1$ and $A_2$ are \emph{trivially equivalent}, if either
$A_1^t=A_2$ with a rational $t=p/q \in\QQ$, with all diagonal
entries $|\lambda_j^{(1)}|\in \NN$ ($j=1,\dots,d$) being full
$q^{\rm th}$ powers (of some, perhaps different natural entries), or
if with some natural numbers $a,b\in\NN$ we have
$|\lambda_j^{(1)}|=a^{n_j}$ with $n_j\in\NN$ ($j=1,\dots,d$)
satisfying $(n_1,\dots,n_d)=q$, and with $t=m/q \cdot \log b/\log
a$, where $m\in \NN$.
\end{defn}

With this notion we can summarize our findings in the next
statement.

\begin{prop}\label{prop:trivieq} Let
$A_1,A_2:\mathbb{R}^{d}\longrightarrow \mathbb{R}^{d}$ be
self-adjoint expansive linear maps, with $A_{\mu}= C_{\mu} J_{\mu}
C_{\mu}^{-1}$, $\mu = 1,2,$ where $J_{\mu}$ are the respective
diagonal maps as in (\ref{Jmu}). Then according to Corollary
\ref{cormain}, $\mathcal{E}_{A_1}=\mathcal{E}_{A_2}$ if and only
if there exists $t>0$ such that $A_1^{'t}=A'_2$, where  $A'_{\mu}=
C_{\mu} J'_{\mu} C_{\mu}^{-1}$, $\mu = 1,2,$ with $J'_{\mu}$ in
(\ref{Jmuprime}). Moreover, if the respective matrices are
\emph{trivially equivalent} in the above sense, then they both
satisfy (\ref{uno}), and thus form two equivalent expansive linear
maps for MRA.
\end{prop}

The next question is to describe solutions of (\ref{uno}) for $a$
and $b$ in the diagonal of a 2 by 2 matrix $A_1$ with linearly
independent logarithms over $\QQ$, and $t\notin\QQ$. We can
conjecture that such equivalences do not occur, i.e. if $A_1$ and
$A_2$ are equivalent positive expansive matrices in $\RR^{2\times
2}$, satisfying (\ref{uno}), then they are from the above
described trivial classes (including, of course, both the cases
when $A_{\mu}$ are dilations, as then $a=b$, and
when $A_1=A_2$, as then $t=1$ is rational).\\

We can not prove this conjecture, but we can say that a well-known
conjecture of number theory would imply this, too. Namely, we can
now recall the so-called ``Four Exponentials Conjecture'', see
e.g. \cite[p.14]{Wa:00}.
\begin{conj}[{Four Exponentials Conjecture}] Let $x_1$, $x_2$
be two $\mathbb{Q}$-linearly independent complex numbers and
$y_1$, $y_2$ also two $\mathbb{Q}$-linearly independent complex
numbers. Then  at least one of the 4 numbers
\[
 \exp(x_i y_j) ,  \qquad (i=1,2, ~~~ j=1,2)
\]
is transcendental.
\end{conj}

Indeed, if the conjecture is right, we can choose $x_1:=\log a$ ,
$x_2:=\log b$, $y_1:= 1$ and $y_2 : = t$. If $x_1$ and $x_2$ are
linearly independent over $\mathbb{Q}$ and $y_1$ and $y_2$ are also
linearly independent over $\QQ$, then either of the four numbers
$a=e^{\log a}$, $b=e^{\log b}$, $a^t=e^{t \log a}$ and $b^t=e^{t
\log b}$ must be transcendental, therefore one can not have
$a,b,a^t,b^t \in \mathbb{Z}$. So in case the Four Exponentials
Conjecture holds true, we must necessarily have either $a=b^q$ where
$q \in \mathbb{Q}$, or $t \in \mathbb{Q}$.

The same argument can be implemented even in dimension $d$.

\begin{prop} Assume that the above Four Exponentials
Conjecture holds true. Then the two self-adjoint expansive linear
maps $A_1$ and $A_2$ generate equivalent MRA if and only if the
conditions of Proposition \ref{prop:trivieq} above hold true:
$A_1^{'t}=A'_2$ and, moreover, $A_1$ and $A_2$ are trivially
equivalent matrices.
\end{prop}

\pf Assume first that $t\notin \QQ$. Then $t$ and $1$ are linearly
independent (over $\QQ$), hence for any pair of indices $1\leq j ,
k\leq d$ applying the Four Exponential Conjecture we conclude
linear dependence of $x_1:=\log|\lambda_j^{(1)}|$ and
$x_2:=\log|\lambda_k^{(1)}|$, i.e.
$|\lambda_k^{(1)}|=|\lambda_j^{(1)}|^{r_{j,k}/s_{j,k}}$, with
$r_{j,k},s_{j,k}\in\NN$. So in view of the unique prime
factorization, $|\lambda_j^{(1)}|=\alpha_{j,k}^{s_{j,k}}$ (and
also $|\lambda_k^{(1)}|=\alpha_{j,k}^{r_{j,k}}$) with some
$\alpha_{j,k}\in\NN$.  Now, for each $k \in \{1,...,d \}$ we
compare the different expressions
$|\lambda_k^{(1)}|=\alpha_{j,k}^{r_{j,k}}$, $j \in \{1,...,d \}$,
thus, taking the l.c.m. of the numbers $r_{j,k}$, $j=1,\dots,d$
we can define $r_k:=[\dots r_{j,k}\dots]$ and, again from the
unique prime factorization, we find that $|\lambda_k^{(1)}|$ must
be a full $r_k$th power, i.e. $|\lambda_k^{(1)}|=\alpha_k^{r_k}$
with some $\alpha_{k}\in\NN$. . Consequently, for any fixed pair
of $k, j \in \{ 1,...,d\}$,
$\alpha_k^{r_k}=|\lambda_k^{(1)}|=|\lambda_j^{(1)}|^{r_{j,k}/s_{j,k}}=
\left(\alpha_j^{r_j}\right)^{r_{j,k}/s_{j,k}}$, so again by the
unique prime factorization all the $\alpha_k$s have the same prime
divisors, and we can write $|\lambda_k^{(1)}|=\alpha^{u_k}$ (with
some $\alpha\in\NN$ and $u_k\in\NN$) for $k=1,\dots,d$. We can
even consider $\nu:=\max \{\mu ~:~ \alpha=a^\mu,~a, \mu\in \NN\}$,
and with this $\nu$ write $|\lambda_k^{(1)}|=a^{n_k}$, where
$n_k:=\nu u_k$ ($k=1,\dots,d$).

Observe that the same reasoning applies to the diagonal entries of
the second matrix $A_2'$, hence we also find
$|\lambda_k^{(2)}|=b^{m_k}$, where $b\in\NN$ and $m_k\in\NN$
($k=1,\dots,d$).

Now write $(n_1,\dots,n_d)=q$, and apply the equivalence condition
$A_1^t=A_2$ to get $t n_k \log a= m_k \log b$ for all $k=1,\dots,d$.
By the linear representation of the g.c.d, we thus obtain $t q \log
a = m \log b$, where $m\in\NN$ is a linear combination of the
exponents $m_k$ (and hence is divisible by $p:=(m_1,\dots,m_d)$). In
any case, we have obtained the case $t=m/q \log a/\log b$ of the
trivial equivalence above.

Second, let $t=p/q\in \QQ$. Then, for each $j=1,\dots,d$, we have
the equation
$|\lambda_j^{(2)}|=|\lambda_j^{(1)}|^t=|\lambda_j^{(1)}|^{p/q}$,
so $|\lambda_j^{(2)}|=a_j^p$ and $|\lambda_j^{(1)}|=a_j^q$, with
$a_j\in\NN$ otherwise arbitrary: and this is the other case of
trivial equivalence, as defined above.

In all, we found that under the assumption of the truth of the
Four Exponentials Conjecture, equivalence with respect to the
notion of $A_\mu$-approximate continuity (or, equivalently,
$A_\mu$-density at 0) and fulfilling condition (\ref{uno}) implies
trivial equivalence of the expansive self-adjoint linear matrices
$A_1$ and $A_2$. Last, but not least, we are indebted to our
referees, who pointed out various connected works and suggested
several further improvements of the originally clumsy
presentation. \qed

\section{Acknowledgements}

We would like to thank Professor Kazaros Kazarian and Patrico
Cifuentes for the useful discussions about our topic, and to
Professor K\'alm\'an Gy\H ory for the information about the Four
Exponentials Conjecture. Also the second author would like to
thank Professor Peter Oswald for allowing the inclusion of the
result of our joint work on the characterization of matrices
equivalent to the dyadic matrix, into his thesis \cite{Thesis}. We
are indebted to two referees, who pointed out several connected
work in the literature and suggested various improvements of the
originally clumsy presentation.


\begin{thebibliography}{9}


\bibitem{B:78}
A. M. Bruckner,  Differentiation of Real Functions,
Springer-Verlag Berlin Heidelberg New York (1978).

\bibitem{CKS:05}
P. Cifuentes, K.S. Kazarian, A.  San Antol\'{\i}n,
Characterization of scaling functions in a multiresolution
analysis,  Proc. Amer. Math. Soc. {\bf 133} No. 4 (2005),
1013-1023.

\bibitem{CKS:06}
P. Cifuentes, K.S. Kazarian, A.  San Antol\'{\i}n,
Characterization of scaling functions,  Wavelets and splines:
Athens 2005, 152--163, Mod. Methods Math., Nashboro Press,
Brentwood, TN, 2006.

\bibitem{DDL:95}
 S. Dahlke, W. Dahmen, V. Latour,  Smooth refinable functions and
 wavelets obtained by convolution
products,
 Appl. Comput. Harmon. Anal. 2 (1995), no. 1,
68--84.

\bibitem{DGL:97} S. Dahlke, K. Gr\"ochenig, V. Latour,
Biorthogonal box spline wavelet bases, in: Surface fitting and
multiresolution methods (Chaminix-Mont Blanc, 1996), 83-92,
Vanderbilt Univ. Press, Nashville, TN, 1997.

\bibitem{DLN:97} S. Dahlke, V. Latour, M. Neeb,
Generalized cardinal $B$-splines: stability, linear independence,
and appropriate scaling matrices, Constr. Approx. 13 (1997), no.
1, 29--56.

\bibitem{DM} S. Dahlke, P. Maass, A note on interpolating scaling
functions, Commun. Appl. Anal. 7 (2003), no. 2-3, 265--279.

\bibitem{DMT} S. Dahlke, P. Maass, G. Teschke, Interpolating
scaling functions with duals, J. Comput. Anal. Appl., 6, (2004),
no. 1, 19-29.

\bibitem{GM:92}
K. Gr\"{o}chening,  W. R. Madych,  Multiresolution analysis, Haar
bases and self-similar tilings of $R^n,$  IEEE Trans. Inform.\
Theory, 38(2) (March 1992), 556--568.

\bibitem{Halmos} P. R. Halmos,  Finite-dimensional vector
spaces, second edition, University Series in Undergraduate
Mathematics, D. Van Nostrand Co. Incl., Princeton, New York.
London, Toronto, 1958.

\bibitem{HW:96}
E. Hern\'{a}ndez, G.  Weiss,  A first course on wavelets, CRC
Press, Inc. 1996.


\bibitem{HK:}
 K. Hoffman, R.  Kunze,  Linear Algebra, Prentice-Hall
Mathematics series, 1961.

\bibitem{LW:95}
J. C. Lagarias, Y. Wang, Haar type orthonormal wavelet bases in
$R\sp 2$,  J. Fourier Anal. Appl. 2 (1995), no. 1, 1--14.

\bibitem{Ma:92}
W. R.  Madych,  Some elementary properties of multiresolution
analyses of $L^2(R^d),$ In { Wavelets - a tutorial in theory and
applications}, Ch. Chui ed., Academic Press (1992), 259--294.

\bibitem{M:89}
S. Mallat,  Multiresolution approximations and wavelet orthonormal
bases of $L^2(R)$, Trans. of Amer. Math. Soc., 315 (1989), 69--87.

\bibitem{Me:90}
 Y. Meyer,  Ondelettes et op$\acute{e}$rateurs. I,
Hermann, Paris (1990) [English Translation: Wavelets and
operators, Cambridge University Press, (1992).]

\bibitem{N:60}
I. P.  Natanson,  Theory of functions of a real variable, London,
vol. I, 1960.

\bibitem{RS:80}
M. Reed,  B. Simon,  Functional Analysis, Academic Press, Inc.
(1980).

\bibitem{Thesis} A. San Antol\'{\i}n,
 Characterization and properties of  scaling functions
and low pass filters of a multiresolution analysis, PhD thesis,
Universidad Autonoma de Madrid, 2007. www.uam.es/angel.sanantolin


\bibitem{S:94}  R. Strichartz,  Construction of orthonormal wavelets,
 Wavelets: mathematics and applications, 23--50,
Stud. Adv. Math., CRC, Boca Raton, FL, 1994.

\bibitem{Wa:00} M. Waldschmidt,  Diophantine approximation on linear algebraic groups.
Transcendence properties of the exponential function in several
variables. Grundlehren der Mathematischen Wissenschaften
[Fundamental Principles of Mathematical Sciences], 326.
Springer-Verlag, Berlin, 2000.

\bibitem{W}
R. Webster,  Convexity, Oxford University Press, Oxford, 1994.

\bibitem{W:97}  P. Wojtaszczyk,  A mathematical
introduction to wavelets, London Mathematical Society, Student
Texts 37 (1997).




\end{thebibliography}
\end{document}